\documentclass[11pt]{article}

\usepackage{graphicx}
\usepackage{amsfonts}

\usepackage{theorem}
\newtheorem{Lem}{Lemma}[section]
\newtheorem{Def}[Lem]{Definition}
\newtheorem{The}[Lem]{Theorem}
\newtheorem{Prop}[Lem]{Proposition}
\newtheorem{Cor}[Lem]{Corollary}

\newtheorem{Rem}[Lem]{Remark}

\newcommand{\qed}{\hbox{\rule{6pt}{6pt}}}
\begin{document}

\title{A note on operator inequalities of Tsallis relative operator entropy}
\author{Shigeru Furuichi$^1$\footnote{E-mail:furuichi@ed.yama.tus.ac.jp}, Kenjiro Yanagi$^2$\footnote{E-mail:yanagi@yamaguchi-u.ac.jp} and Ken Kuriyama$^2$\footnote{E-mail:kuriyama@yamaguchi-u.ac.jp}\\
$^1${\small Department of Electronics and Computer Science,}\\{\small Tokyo University of Science, Onoda City, Yamaguchi, 756-0884, Japan}\\
$^2${\small Department of Applied Science, Faculty of Engineering,} \\{\small 
Yamaguchi University,Tokiwadai 2-16-1, Ube City, 755-0811, Japan}}
\date{}
\maketitle

{\bf Abstract.}
Tsallis relative operator entropy was defined as a parametric extension of relative operator entropy and the generalized Shannon inequalities were shown in the previous paper.
After the review of some fundamental properties of Tsallis relative operator entropy, some operator inequalities related to Tsallis relative operator
entropy are shown in the present paper. Our inequalities give the upper and lower bounds of Tsallis relative operator entropy. 
The operator equality on Tsallis relative operator entropy is also shown by considering the tensor product. This relation generalizes the pseudoadditivity
for Tsallis entropy. As a corollary of our operator equality derived from the tensor product manipulation, we show several operator inequalities including the superadditivity and the subadditivity for
Tsallis relative operator entropy. Our results are generalizations of the superadditivity and the subadditivity for Tsallis entropy.
\vspace{3mm}

{\bf Keywords : } Tsallis relative entropy, relative operator entropy, superadditivity, subadditivity and operator inequality

\vspace{3mm}

{\bf 2000 Mathematics Subject Classification : } 47A63, 94A17, 15A39

\section{Introduction}

The relative operator entropy 
\begin{equation}\label{roe}
S(A\vert B) \equiv A^{1/2}\log(A^{-1/2}BA^{-1/2})A^{1/2}
\end{equation}
for two invertible positive operators $A$ and $B$,
was introduced by Fujii and Kamei in \cite{FuKa:rel}, as a generalization of the operator entropy 
\begin{equation} \label{umegaki_nakamura} 
H(A) \equiv S(A\vert I)= -A\log A
\end{equation}
defined in \cite{NU}.
In the present paper, we study a parametric extension of relative operator entropy, we call it Tsallis relative operator entropy which was firstly
introduced in our previous paper \cite{YKF}, in the following manner.
\begin{Def} {\bf (\cite{YKF})}
For two invertible positive operators $A$ and $B$ on Hilbert space, and any real number $\lambda \in (0,1]$, Tsallis relative operator entropy is
defined by
\begin{equation}\label{troe}
T_{\lambda}(A\vert B) \equiv  \frac{A^{1/2} (A^{-1/2}BA^{-1/2})^{\lambda} A^{1/2} -A}{\lambda}.
\end{equation}
\end{Def}

Note that more general operator than the above Eq.(\ref{troe}) has been introduced in \cite{FK}, it was named {\it solidarity} and then several properties were shown in \cite{FK,FFS,F}.

In statistical physics, C.Tsallis \cite{Ts:rel} introduced the parametrically extended Shannon entropy to study the multifractal system. 
It is called Tsallis entropy and is defined for the probability distribution $p(x)\equiv p(X=x)$ of the random variable $X$ such as
$$
S_{q}(X)\equiv \frac{1- \sum_x p(x)^q}{q-1}, \quad ( q \in \mathbf{R}).
$$
The definition of Tsallis entropy was generalized for density operator $\rho$ (positive operator with unit trace) in noncommutative (quantum) case, such that 
$$
S_{\lambda}   (\rho) \equiv \frac{Tr[\rho^{1-\lambda}]-1}{\lambda},\quad (0 < \lambda \leq 1).
$$
We use the parameter $\lambda$ in the present paper, instead of the parameter $q$ which is often used to study Tsallis entropy, where we have the simple relation $\lambda = 1-q$. 
As main features of Tsallis entropy, it has 
\begin{itemize}
\item[(i)] the concavity in $\rho$ for all $0 < \lambda \leq 1$,
\item[(ii)] the pseudoadditivity 
\begin{equation}  \label{pseudo_Tsallis}
S_{\lambda}(\rho_1\otimes \rho_2) =S_{\lambda}(\rho_1)+S_{\lambda}(\rho_2)+\lambda S_{ \lambda}(\rho_1)S_{\lambda}(\rho_2),
\end{equation}
and 
\item[(iii)] it converges to von Neumann entropy $S(\rho) \equiv Tr[H(\rho)]$ as $ \lambda \to 0$. 
\end{itemize}
It is remarkable that R\'enyi entropy (which is also called $\alpha$ entropy) \cite{Re} 
$\frac{1}{1-\alpha}\log \sum_x p(x)^{\alpha}$ has also the concavity in $p(x)$ for the limited region of the parameter $\alpha < 1$,
but it lacks the concavity in the case of $\alpha > 1$, which means that it is not suitable to use for the study of thermodynamical physics.
(The classical Tsallis entropy $S_q(X)$ has the concavity in $p(x)$ for all $q \geq 0$.)
There is a different mathematical property between R\'enyi entropy and Tsallis entropy in the sense that Tsallis entropy has a pseudoadditivity
while R\'enyi entropy has an additivity.
After the birth of Tsallis entropy, the {\it nonadditive} (sometimes called {\it nonextensive} or {\it Tsallis}) statistical physics has been in 
progresses with Tsallis entropy and the related entropic quantities such as Tsallis relative entropy \cite{AO:re1,Ab:re1,FYK}.
Therefore it is important to study the properties of such entropic quantitties in mathematical physics. 
We studied the fundamental properties of Tsallis relative entropy from the information theoretical point of view, 
and derived some important trace inequalities in \cite{FYK}.
The Tsallis relative entropy in quantum system was defined by
\begin{equation}
D_{ \lambda} (\rho\vert \sigma) \equiv \frac{1-Tr[\rho^{1-\lambda}\sigma^{\lambda}]}{\lambda}
\end{equation}
for any $0 < \lambda \leq 1$ and  any two density operators $\rho$ and $\sigma$.
It has similar features to Tsallis entropy such that
 \begin{itemize}
\item[(i)] $D_{ \lambda}(\rho\vert \sigma) $ has the joint convexity in $\rho$ and $\sigma$ and for all $0 < \lambda \leq 1$, 
\item[(ii)] $D_{ \lambda}(\rho\vert \sigma) $ has the pseudoadditivity 
\begin{equation}  \label{pseudo_relative_Tsallis}
D_{ \lambda}  (\rho_1 \otimes \rho_2 \vert \sigma_1 \otimes \sigma_2) =D_{ \lambda}  (\rho_1\vert \sigma_1)+D_{ \lambda}  (\rho_2\vert \sigma_2 )
-\lambda D_{ \lambda}  (\rho_1\vert \sigma_1)D_{ \lambda}  (\rho_2\vert \sigma_2),
\end{equation}
and 
\item[(iii)] it converges to the quantum relative entropy $U(\rho \vert \sigma) \equiv Tr[\rho (\log \rho -\log \sigma)]$ introduced by Umegaki in \cite{Ume},
 by taking the limit as $\lambda \to 0$. 
\end{itemize}
As it was shown in our previour paper \cite{FYK},  we have the trace inequality,
\begin{equation} \label{gHP}
D_{ \lambda}(\rho \vert \sigma) \leq -Tr[T_{\lambda}(\rho \vert \sigma)] 
\end{equation}
which means that the minus of the trace of Tsallis relative operator entropy for two density operators $\rho$ and $\sigma$ is not equal to
Tsallis relative entropy in general. This trace inequality is a generalization of $U(\rho\vert \sigma) \leq -Tr[S(\rho\vert \sigma)]$ proved by
Hiai and Petz in \cite{HP}. If $\rho$ and $\sigma$ are commutative, then we have 
\begin{equation} \label{gHPeq}
D_{\lambda}(\rho\vert \sigma) = -Tr[T_{\lambda}(\rho\vert \sigma)].
\end{equation}

We are now interested in the operator inequalities \cite{An,Furuta,Zh}, before taking
the trace of Tsallis relative entropy.
As it is meaningful to study relative operator entropy for the development of the noncommutative information theory,
we think it is indispensable to study the properties of Tsallis relative operator entropy for the progress of noncommutative
statisitical physics and nonadditive quantum information theory \cite{AO:re1}. 

From the definition, we easily find that $\lim_{\lambda \to 0}T_{\lambda}(A\vert B) =S(A\vert B)$, 
since $\ln_{\lambda}(x)\equiv \frac{x^{\lambda}-1}{\lambda}$ uniformly converges to
$\log x$ for $x \geq 0$ by Dini's theorem. Also we have $H_{\lambda}(A)\equiv  T_{\lambda}(A\vert I) =  \frac{A^{1-\lambda}-A}{\lambda}$, 
which we call it Tsallis operator entropy. We easily find $H_{\lambda}(A)$ 
converges to the operator entropy $H(A)$ as $\lambda \to 0$.
In the following section, we study the basical properties of Tsallis relative operator entropy.

\section{Review of properties of $T_{\lambda}(A\vert B)$}

In this section, we give several fundamental properties of Tsallis relative operator entropy along the line of the paper \cite{FuKa:rel}.
More general results were shown in \cite{FK,FFS,F}. However, in this short paper, 
for the convenience of the readers, we will review some properties of Tsallis relative operator
entropy which is a special case of the solidarity introduced in \cite{FK}.
As a matter of convenience, we often use the notation $\ln_{\lambda}X \equiv \frac{X^{\lambda}-I}{\lambda}$ for the positive operator $X$,
 through this paper.
Then we often rewrite Tsallis relative operator entropy $T_{\lambda}(A \vert B)$ as
$$ T_{\lambda}(A \vert B) =A^{1/2}\ln_{\lambda}(A^{-1/2}BA^{-1/2})A^{1/2}. $$
All results in this section recover the properties of relative operator entropy in \cite{FuKa:rel}  as $\lambda \to 0$.
The Tsallis relative operator entropy $T_{\lambda}(A \vert B) $ has the following properties.
\begin{Prop}\label{a1}{\bf (\cite{FFS})}
\begin{itemize}
\item[(1)] (homogeneity) $T_{\lambda}(\alpha A \vert \alpha B) =\alpha T_{\lambda}(A\vert B)$  for any positive number $\alpha$.
\item[(2)] (monotonicity) If $B \leq C$, then $T_{\lambda}(A\vert B) \leq T_{\lambda}(A\vert C).$
\end{itemize}
\end{Prop}
{\bf (Proof)}
(1) is trivial.
Through this paper, we use the following notation 
$$A \sharp_{\lambda} B \equiv A^{1/2} (A^{-1/2}BA^{-1/2})^{\lambda} A^{1/2}$$
representing the operator mean between $A$ and $B$, which is often called $\lambda$-power mean.
See \cite{HK} for the operator mean.
By the use of this notation, Tsallis relative operator entropy can be rewritten by 
$T_{\lambda}(A \vert B) = \frac{1}{\lambda}(A \sharp_{\lambda}  B -A)$.
It is known that the monotonicity on the operator mean such that if $A\leq B$ and $C \leq D$, 
then $A \sharp_{\lambda}  B \leq C \sharp_{\lambda}  D$.
Therefore if $B \leq C$, then $T_{\lambda}(A\vert B) \leq T_{\lambda}(A\vert C).$
\hfill \qed

\vspace{3mm}

Note that the general results were shown in \cite{FFS}.
Moreover, $T_{\lambda}(A\vert B)$ has the superadditivity and joint concavity.

\begin{Prop}
\begin{itemize}
\item[(1)] (superadditivity) $T_{\lambda}(A_1+A_2\vert B_1+B_2) \geq T_{\lambda}(A_1\vert B_1) + T_{\lambda}(A_2 \vert B_2).$
\item[(2)] (joint concavity) $T_{\lambda}(\alpha A_1+\beta A_2\vert \alpha B_1+\beta B_2) \geq \alpha T_{\lambda}(A_1\vert B_1) + \beta T_{\lambda}(A_2 \vert B_2) .$
\end{itemize}
\end{Prop}
{\bf (Proof)}
As mentioned above, we can rewrite $T_{\lambda}(A\vert B) =\frac{1}{\lambda}(A \sharp_{\lambda}  B -A) $ by means of the operator mean 
$A \sharp_{\lambda}  B$.
It is known that in general the concavity $(A+B) m  (C+D) \geq A m C +B m D$ holds for the operator mean $m$.
Therefore we have 
\begin{eqnarray*}
T_{\lambda}(A_1+A_2\vert B_1+B_2) &=& \frac{1}{\lambda} \left\{ (A_1 +A_2) \sharp_{\lambda}  (B_1 +B_2) -(A_1+A_2) \right\} \\
& \geq & \frac{1}{\lambda} \left( A_1 \sharp_{\lambda}  B_1 - A_1 +A_2 \sharp_{\lambda}  B_2 -A_2 \right) \\
&=& T_{\lambda}(A_1\vert B_1) + T_{\lambda}(A_2 \vert B_2).
 \end{eqnarray*}
(2) follows from the superadditivity and homogeneity.
\hfill \qed

\vspace{3mm}

We also easily find that we have
$$T_{\lambda}(UAU^*\vert UBU^*) = U T_{\lambda} (A\vert B) U^*$$ 
for any unitary operator $U$.

Finally we show the monotonicity of Tsallis relative operator entropy.

\begin{Prop} {\bf (\cite{F})}
For a unital positive linear map $\Phi$ from the set of the bounded linear operators on Hilbert space to itself, we have
$$\Phi (T_{\lambda}(A\vert B)) \leq T_{\lambda} (\Phi (A) \vert \Phi (B)) .$$
\end{Prop}

{\bf (Proof)}
Putting $X=A^{1/2}B^{-1/2}$ and the formula $\ln_{\lambda}t =-t^{\lambda} \ln_{\lambda} t^{-1}$, Tsallis relative operator entropy
is calculated as
\begin{eqnarray*}
T_{\lambda}(A\vert B) &=& A^{1/2} (\ln_{\lambda} A^{-1/2} B A^{-1/2}   ) A^{1/2} \\
&=& A^{1/2}\left\{   -(A^{-1/2} B A^{-1/2})^{\lambda}\ln_{\lambda} (A^{1/2} B^{-1} A^{1/2})  \right\} A^{1/2} \\
&=& A^{1/2}\left\{   -(XX^*)^{-\lambda}\ln_{\lambda} (XX^*)  \right\} A^{1/2} \\
&=& A^{1/2}\left\{   -(XX^*)^{-\lambda}\ln_{\lambda} (XX^*)X  \right\}  B^{1/2} 
\end{eqnarray*}

In general, for any bounded linear operator $A$ on Hilbert space and the continuous function $f$ on the interval
 $[0, \vert\vert A\vert\vert^2]$,
we have $Af(A^*A)= f(AA^*)A$. Thus we have
\begin{eqnarray*}
T_{\lambda}(A\vert B) &=&A^{1/2}\left\{   -X (X^*X)^{-\lambda}\ln_{\lambda} (X^*X)  \right\}  B^{1/2} \\
&=& B^{1/2} \left\{ -X^*X (X^*X)^{-\lambda}\ln_{\lambda} (X^*X)  \right\}  B^{1/2} \\
&=& B^{1/2} \left\{ F(X^*X)\right\}  B^{1/2} \\
&=& B^{1/2} \left\{ F(B^{-1/2} AB^{-1/2})\right\}  B^{1/2}, 
\end{eqnarray*}
where $F$ is the concave function defined by $F(t) =-t^{1-\lambda} \ln_{\lambda} t , (0 < \lambda \leq 1)$.
Therefore Tsallis relative operator entropy can be written by 
$$ T_{\lambda}(A\vert B) = F(A \slash B ) \slash B^{-1},$$
where $X\slash Y \equiv Y^{-1/2}XY^{-1/2}$, which is same notation used in \cite{FuKa:rel}. 
Thus the present theorem follows by the similar way of the proof of Theorem 7 in \cite{FuKa:rel}.
\hfill \qed

\vspace{3mm}

Note that the general results were shown in \cite{F}.

\section{Inequalities as bounds for $T_{\lambda}(A\vert B)$}

In the previous section, some basical properties of Tsallis relative operator entropy are shown.
In this section, we discuss on the upper and lower bounds of Tsallis relative operator entropy.

Firstly we prove the relation of ordering between Tsallis relative operator entropy and relative operator entropy.

\begin{Prop} {\bf (\cite{FK})}
For any invertible positive operator $A$ and $B$, $0< \lambda \leq 1$, we have
$$
T_{-\lambda}(A\vert B) \leq S(A\vert B) \leq T_{\lambda}(A\vert B).
$$
\end{Prop}

{\bf (Proof)}
Since we have 
$$
\frac{x^{-\lambda}-1}{-\lambda} \leq \log x \leq \frac{x^{\lambda}-1}{\lambda},
$$
for any $x >0 $ and $\lambda >0$,
we have the following inequalities
$$
\ln_{-\lambda}(A^{-1/2}BA^{-1/2})  \leq \log(A^{-1/2}BA^{-1/2})  \leq \ln_{\lambda}(A^{-1/2}BA^{-1/2}).
$$
Multiplying $A^{1/2}$ to both sides from both sides, the present proposition follows.
\hfill \qed

\vspace{3mm}

As for the bounds of Tsallis relative operator entropy, we have the following results.
\begin{Prop}
\begin{itemize}
\item[(1)] $T_{\lambda}(A\vert B) \leq H_{\lambda}(A) + A^{1-\lambda} \ln_{\lambda} \vert\vert B\vert\vert .$
\item[(2)] If $\mu A\leq B$, then $T_{\lambda} (A \vert B) \geq (\ln_{\lambda} \mu ) A.$
\end{itemize}
\end{Prop}
{\bf (Proof)}
(1) follows from the direct calculations, since $B \leq \vert \vert B\vert \vert I$.
(2) also follows from the direct calculations, since $\mu A \leq B.$ 
\hfill \qed

\vspace{3mm}

Moreover we have the following bounds for Tsallis relative operator entropy.

\begin{Lem}   \label{3lem1}
For any positive real number $x$ and $0< \lambda \leq 1$, the following inequalities hold
$$
1-\frac{1}{x} \leq \ln_{\lambda}x \leq x-1.
$$
\end{Lem}
{\bf (Proof)}
It follows from easy calculations.
\hfill \qed

\vspace{3mm}

\begin{Prop} {\bf (\cite{FK})}
For any invertible positive operator $A$ and $B$, $0< \lambda \leq 1$,
\begin{equation} \label{the3-4}
A-AB^{-1}A \leq T_{\lambda}(A\vert B) \leq -A+B. 
\end{equation}
Moreover, $T_{\lambda} (A\vert B) =0$ if and only if $A=B$.
\end{Prop}
{\bf (Proof)}
From Lemma \ref{3lem1}, we have
$$
I-A^{1/2}B^{-1}A^{1/2} \leq \ln_{\lambda}(A^{-1/2}BA^{-1/2}) \leq -I +A^{-1/2}BA^{-1/2}.
$$
Multiplying $A^{1/2}$ from the both sides, we have
$$
A^{1/2} (I-A^{1/2}B^{-1}A^{1/2})A^{1/2} \leq A^{1/2} \ln_{\lambda}(A^{-1/2}BA^{-1/2}) A^{1/2} \leq A^{1/2} (-I +A^{-1/2}BA^{-1/2})A^{1/2} .
$$
Therefore we have
$$A-AB^{-1}A \leq T_{\lambda}(A\vert B) \leq -A+B. $$

Moreover, we suppose $T_{\lambda}(A\vert B)=0$. Then, from the above inequalities, we have
$$A-AB^{-1}A \leq 0 \leq B-A,$$
which implies $A\leq B$ and $A \geq B$. Thus we have $A=B$.
If $A=B$, then we easily find $T_{\lambda}(A\vert B)=0.$
\hfill \qed

\vspace{3mm}

Finally we prove the further bounds of Tsallis relative operator entropy with the $\lambda$-power mean $\sharp_{\lambda}$.
\begin{Lem}\label{finallemma}
For any positive real number $\alpha$ and $x$ and $0< \lambda \leq 1$, the inequalities hold
\begin{equation}\label{finlem}
x^{\lambda}\left(1-\frac{1}{\alpha x}\right) + \ln_{\lambda}\frac{1}{\alpha} \leq \ln_{\lambda} x \leq \frac{x}{\alpha} -1 -x^{\lambda}\ln_{\lambda}\frac{1}{\alpha}.
\end{equation}
The equality of the right hand side of the above inequalities hold if and only if $x=\alpha$.
The equality of the left hand side of the above inequalities hold if and only if $x=\frac{1}{\alpha}$.
\end{Lem}
{\bf (Proof)}
Since we have $\ln_{\lambda} z \leq z-1$ for any $z>0$ and $0< \lambda\leq 1$,
we have $\ln_{\lambda}\frac{x}{\alpha} \leq \frac{x}{\alpha}-1$ for any $\alpha >0 $ and $x>0$. 
By the formula $\ln_{\lambda}\frac{x}{y} = \ln_{\lambda}x+x^{\lambda}\ln_{\lambda}\frac{1}{y}$,
we have 
$$\ln_{\lambda}x \leq \frac{x}{\alpha} -1 -x^{\lambda}\ln_{\lambda}\frac{1}{\alpha},$$
which implies the right hand side of the inequalities Eq.(\ref{finlem}).
Putting $\frac{1}{x}$ instead of $x$ in the above inequality, we have the left hand of inequalities Eq.(\ref{finlem}), 
since $\ln_{\lambda}\frac{1}{x}=-x^{-\lambda}\ln_{\lambda}x$.
\hfill \qed

\vspace{3mm}

\begin{The}
For any invertible positive operators $A$ and $B$, and any positive real number $\alpha$, the following inequality holds
\begin{equation}\label{alphax}
A\sharp_{\lambda}B -\frac{1}{\alpha} A\sharp_{\lambda -1}B + (\ln_{\lambda} \frac{1}{\alpha}) A
\leq T_{\lambda}(A\vert B) \leq \frac{1}{\alpha} B -A - (\ln_{\lambda}\frac{1}{\alpha}) A\sharp_{\lambda}B.
\end{equation}
The equality of the right hand side of the above inequalities holds if and only if $B = \alpha A$.
The equality of the left hand side of the above inequalities holds if and only if $A = \alpha B$.
We have that $T_{\lambda}(A\vert B)=0$ is equivalent to $A=B$.
\end{The}
{\bf (Proof)}
From Lemma \ref{finallemma}, we have
\begin{eqnarray*}
& & (A^{-1/2} B A^{-1/2})^{\lambda} - \frac{1}{\alpha} (A^{-1/2} B A^{-1/2})^{\lambda -1} + (\ln_{\lambda}\frac{1}{\alpha}) I \leq
\ln_{\lambda}(A^{-1/2} B A^{-1/2}) \\
& &\leq \frac{1}{\alpha} (A^{-1/2} B A^{-1/2}) -I - (\ln_{\lambda}\frac{1}{\alpha}) (A^{-1/2} B A^{-1/2})^{\lambda}.
\end{eqnarray*}
Multiplying $A^{1/2}$ to both sides from the both sides, we have the present theorem. 
It is clear that the equality conditions follow from Lemma \ref{finallemma}.

We assume that $T_{\lambda}(A\vert B)=0$ holds. Then we have 
$$A\sharp_{\lambda}B -\frac{1}{\alpha} A\sharp_{\lambda -1}B + (\ln_{\lambda} \frac{1}{\alpha}) A
\leq 0 \leq \frac{1}{\alpha} B -A - (\ln_{\lambda}\frac{1}{\alpha}) A\sharp_{\lambda}B.$$
Putting $\alpha = 1$ in the right hand side of the avobe inequalities, we have $A \leq B$. Moreove
putting $\alpha =1$ and $\lambda =1 $ in the left hand side of the above inequalities, we have $B \leq A$. 
Thus we have $A=B$. Conversely, if $A=B$ then $T_{\lambda}(A\vert B)=0$ easily follows from the definition of $T_{\lambda}(A\vert B)$. 
\hfill \qed

\vspace{3mm}

\begin{Rem}
We note that Eq.(\ref{alphax}) recovers the inequalities shown in \cite{Furuta} :
$$
(1-\log \alpha)A -\frac{1}{\alpha}AB^{-1}A \leq S(A\vert B) \leq (\log \alpha -1)A + \frac{1}{\alpha} B
$$
as $\lambda \to 0$.
Moreover, if we put $\alpha = 1$, then we have 
$$
A -AB^{-1}A \leq S(A\vert B) \leq B-A
$$
which recover the inequalities of Corollary 5 in \cite{FK}, cf. Eq.(\ref{the3-4}).
\end{Rem}

\section{Inequalities on $T_{\lambda}(A\vert B)$ derived from the tensor product}

Taking account for the pseudoadditivity (Eq.(\ref{pseudo_Tsallis}) and Eq.(\ref{pseudo_relative_Tsallis})) which are the typical features of Tsallis entropies, we consider the
Tsallis relative operator entropy of two positive operator of the tensor product $A_1 \otimes A_2$ and $B_1 \otimes B_2$. 
To show our theorem, we state the following lemma for the convenience of the readers.

\begin{Lem}  \label{lem_tensor}
For any real number $a$ and any strictly positive operators $X$ and $Y$, we have  $$(X\otimes Y)^{a} = X^a \otimes Y^a. $$ 
\end{Lem}
{\bf (Proof)}
Firstly we have 
\begin{equation}   \label{eq_lem_ten1}
(X\otimes Y)^{n} = X^n \otimes Y^n
\end{equation}
 for any natural number $n$ and any strictly positive operators $X$ and $Y$.
Since we have $( X^{-1} \otimes Y^{-1} ) (X \otimes Y ) =I \otimes I$, we also have 
\begin{equation}   \label{eq_lem_ten2}
(X\otimes Y)^{-1} = X^{-1} \otimes Y^{-1}.
\end{equation}
By Eq.(\ref{eq_lem_ten1}), we have $W\otimes Z =\left( W^m \otimes Z^m\right)^{1/m}$ for any natural number $m$ and any strictly  positive operators $W$ and $Z$. 
Putting $W=X^{1/m}$ and $Z=Y^{1/m}$, we have 
\begin{equation}   \label{eq_lem_ten3}
(X\otimes Y)^{1/m}=   X^{1/m} \otimes Y^{1/m}.
\end{equation}
From Eq.(\ref{eq_lem_ten1}), Eq.(\ref{eq_lem_ten2}) and Eq.(\ref{eq_lem_ten3}),  we thus have $(X\otimes Y)^q = X^q \otimes Y^q$ for any
rational number $q$ and any strictly positive operators $X$ and $Y$. 
Therefore we have the present lemma by the fact that any real number can be approximated by a rational munber.
\hfill \qed

\vspace{3mm}

\begin{The}   \label{the_tensor}
For any $0 < \lambda \leq 1$ and any strictly positive operators $A_1,A_2,B_1$ and $B_2$, we have
\begin{equation}  \label{eq_tensor}
T_{\lambda}(A_1\otimes A_2\vert B_1\otimes B_2)=  T_{\lambda}(A_1\vert B_1)\otimes A_2 +A_1 \otimes T_{\lambda}(A_2\vert B_2)
+\lambda T_{\lambda}(A_1\vert B_1) \otimes T_{\lambda}(A_2\vert B_2).
\end{equation}
\end{The}
{\bf (Proof)}
From Lemma \ref{lem_tensor}, we have $(A_1\otimes A_2) \sharp_{\lambda} (B_1 \otimes B_2) = \left( A_1 \sharp_{\lambda} B_1 \right) \otimes \left( A_2 \sharp_{\lambda} B_2 \right)$. 
Then we directly calculate
\begin{eqnarray*}
&& T_{\lambda}(A_1\otimes A_2 \vert B_1 \otimes B_2) = \frac{1}{\lambda} \left\{ \left(A_1\otimes A_2\right) \sharp_{\lambda} \left(B_1 \otimes B_2\right) - A_1 \otimes A_2 \right\} \\
&& = \frac{1}{\lambda} \left\{ \left( A_1 \sharp_{\lambda} B_1 \right) \otimes \left( A_2 \sharp_{\lambda} B_2 \right) - A_1 \otimes A_2 \right\} \\
&& = \frac{1}{\lambda} \left\{  \frac{1}{2} \left( A_1 \sharp_{\lambda} B_1 \right) \otimes \left( A_2 \sharp_{\lambda} B_2 \right) - \frac{1}{2} A_1 \otimes \left( A_2 \sharp_{\lambda} B_2 \right) 
+ \frac{1}{2} \left( A_1 \sharp_{\lambda} B_1 \right) \otimes \left( A_2 \sharp_{\lambda} B_2 \right) \right. \\
&& \left. - \frac{1}{2} \left( A_1 \sharp_{\lambda} B_1 \right)  \otimes A_2 + \frac{1}{2} A_1 \otimes \left( A_2 \sharp_{\lambda} B_2 \right) - \frac{1}{2} A_1\otimes A_2 +\frac{1}{2} \left( A_1 \sharp_{\lambda} B_1 \right)  \otimes A_2  \right. \\
&& \left.  - \frac{1}{2} A_1\otimes A_2  \right\}\\
&& = \frac{1}{2} T_{\lambda} (A_1 \vert B_1) \otimes \left( A_2 \sharp_{\lambda} B_2 \right) +\frac{1}{2}  \left( A_1 \sharp_{\lambda} B_1 \right) \otimes  T_{\lambda} (A_2 \vert B_2) 
+  \frac{1}{2} A_1 \otimes T_{\lambda} (A_2 \vert B_2) \\
&& + \frac{1}{2}  T_{\lambda} (A_1 \vert B_1) \otimes  A_2 \\
&& = \frac{1}{2} T_{\lambda} (A_1 \vert B_1) \otimes \left( A_2 \sharp_{\lambda} B_2 \right) - \frac{1}{2}  T_{\lambda} (A_1 \vert B_1) \otimes  A_2 +T_{\lambda} (A_1 \vert B_1) \otimes  A_2 \\
&& + \frac{1}{2} \left( A_1 \sharp_{\lambda} B_1 \right) \otimes  T_{\lambda} (A_2 \vert B_2) -  \frac{1}{2} A_1 \otimes T_{\lambda} (A_2 \vert B_2) + A_1 \otimes T_{\lambda} (A_2 \vert B_2) \\
&& = \frac{\lambda}{2} T_{\lambda}(A_1\vert B_1) \otimes T_{\lambda}(A_2\vert B_2) +\frac{\lambda}{2} T_{\lambda}(A_1\vert B_1) \otimes T_{\lambda}(A_2\vert B_2) + T_{\lambda}(A_1\vert B_1)\otimes A_2 \\
&& +A_1 \otimes T_{\lambda}(A_2\vert B_2) \\
&&= \lambda T_{\lambda}(A_1\vert B_1) \otimes T_{\lambda}(A_2\vert B_2) + T_{\lambda}(A_1\vert B_1)\otimes A_2 +A_1 \otimes T_{\lambda}(A_2\vert B_2).
\end{eqnarray*}
\hfill \qed

\vspace{3mm}

\begin{Rem}
Theorem \ref{the_tensor} can be also proven by the use of the equality:
$$\ln_{\lambda}(X\otimes Y) = (\ln_{\lambda}X)\otimes I +I\otimes (\ln_{\lambda}Y) +\lambda (\ln_{\lambda}X) \otimes(\ln_{\lambda}Y),$$
for any $0 < \lambda \leq 1$ and any positive operators $X$ and $Y$.
\end{Rem}

Taking the limit as $\lambda \to 0$ in Theorem \ref{the_tensor}, we have the following corollary.
\begin{Cor} \label{cor_tensor}
For any strictly positive operators $A_1,A_2,B_1$ and $B_2$, we have
$$
S(A_1 \otimes A_2 \vert B_1 \otimes B_2) = S(A_1 \vert B_1) \otimes A_2 + A_1 \otimes S(A_2 \vert B_2).
$$
\end{Cor} 

\begin{Rem}
Corollary \ref{cor_tensor} can be also proven by the use of  of the equality:
$$\log (X\otimes Y) = (\log X) \otimes I + I \otimes (\log Y),$$
for any positive operators $X$ and $Y$.
\end{Rem}

\begin{Cor} 
For any $0 < \lambda \leq 1$ and any invertible density operators $\rho_1$ and $\rho_2$, we have Eq.(\ref{pseudo_Tsallis}).
\end{Cor}
{\bf (Proof)}
We put $B_1=B_2=I$ in Theorem \ref{the_tensor}. Then we have
$$H_{\lambda}(A_1 \otimes A_2) = H_{\lambda}(A_1) \otimes A_2 + A_1 \otimes H_{\lambda}(A_2)  + \lambda H_{\lambda}(A_1) \otimes H_{\lambda}(A_2). $$
We put $A_i = \rho_i, (i=1,2) $ in the above and take the trace on the both side, then we have our claim, since $S_{\lambda}(\rho) = Tr[H_{\lambda}(\rho)]$ for density operator $\rho$.
\hfill \qed

\vspace{3mm}

\begin{Rem}
The additivity of von Neumann entropy  $S(\rho)$ :
$$
S(\rho_1 \otimes \rho_2) = S(\rho_1) + S(\rho_2)
$$
is recovered in the limit of Eq.(\ref{pseudo_Tsallis}) as $\lambda \to 0$. Since the quantum Tsallis entropy is nonnegative for any $0 < \lambda \leq 1$ and any density operators $\rho_i, (i=1,2)$, we have
the superadditivity for the quantum Tsallis entropy :
\begin{equation}  \label{cor2}
S_{\lambda}(\rho_1 \otimes \rho_2) \geq S_{\lambda}(\rho_1) + S_{\lambda}( \rho_2)
\end{equation}
holds.
We note that the pseudoadditivity Eq.(\ref{pseudo_Tsallis}) is a special case of 
the pseudoadditivity for the quantum Tsallis relative entropy Eq.(\ref{pseudo_relative_Tsallis}) . Theorem \ref{the_tensor} does not directly imply the pseudoadditivity for the quantum
relative entropy, due to inequality Eq.(\ref{gHP}).
However, in commutative (classical) case, Eq.(\ref{eq_tensor}) directly recovers the pseudoadditivity (Proposition 1.2 (4) in \cite{FYK}) of Tsallis relative entropy for
probability distributions, due to Eq.(\ref{gHPeq}). 
(See \cite{FYK} on the relation between the quantum Tsallis relative entropy and Tsallis relative operator entropy.)
\end{Rem}

Since we have $T_{\lambda}(A\vert B) \geq 0$ for any $B \geq A$, $T_{\lambda}(A\vert B) \leq 0$ for any $B \leq A$ and we have
$X\otimes Y \geq 0$ for any $X \geq 0$ and $Y \geq 0$, Theorem \ref{the_tensor} implies the following corollary.
\begin{Cor}   \label{fininal_cor}
\begin{itemize}
\item[(1)] For any $0 < \lambda \leq 1$ and $0< A_i \leq B_i$, $(i=1,2)$, we have the following inequalities.
\begin{itemize}
\item[(a)] $T_{\lambda}(A_1\otimes A_2 \vert B_1 \otimes B_2) \geq \lambda T_{\lambda}(A_1 \vert B_1)\otimes T_{\lambda}( A_2 \vert B_2).  $
\item[(b)] $T_{\lambda}(A_1\otimes A_2 \vert B_1 \otimes B_2) \geq  T_{\lambda}(A_1 \vert B_1)\otimes A_2 + A_1 \otimes T_{\lambda}( A_2 \vert B_2).$
\end{itemize}
\item[(2)]  For any $0 < \lambda \leq 1$ and $0< B_i \leq A_i$, $(i=1,2)$, we have the following inequalities.
\begin{itemize}
\item[(c)] $T_{\lambda}(A_1\otimes A_2 \vert B_1 \otimes B_2) \leq \lambda T_{\lambda}(A_1 \vert B_1)\otimes T_{\lambda}( A_2 \vert B_2).  $
\item[(d)] $T_{\lambda}(A_1\otimes A_2 \vert B_1 \otimes B_2) \geq  T_{\lambda}(A_1 \vert B_1)\otimes A_2 + A_1 \otimes T_{\lambda}( A_2 \vert B_2).$
\end{itemize}
\end{itemize}
\end{Cor}

\begin{Rem}  \label{final_rem}
We easily find from the above (b) and (d) that we have the superadditivity without depending on the ordering $A_i$  and $B_i$. 
This superadditivity for Tsallis relative operator entropy is a generalization of the superadditivity for Tsallis entropy Eq.(\ref{cor2}). 
Indeed, the inequality Eq.(\ref{cor2}) follows if we put $B_1=B_2=I$ and $A_i= \rho_i,(i=1,2) $, where $\rho_i$ are invertible density operators, 
 in  the above (b) and (d) and then take the trace of them.
\end{Rem}

\section{Concluding remarks}

If we define Tsallis relative operator entropy for the parameter $\lambda \leq 0$ and two invertible positive operators $A$ and $B$ by
$$\widetilde{T_{\lambda}}(A\vert B) \equiv \frac{A^{1/2}(A^{-1/2}BA^{-1/2})^{\lambda}A^{1/2} -A}{\lambda},$$
we have the following operator inequalities by the similar way of Corollary \ref{fininal_cor}.

\begin{Prop}   \label{final_prop}
\begin{itemize}
\item[(1)] For any $\lambda \leq 0$ and $0< A_i \leq B_i$, $(i=1,2)$, we have the following inequalities.
\begin{itemize}
\item[(a')] $  \widetilde{T_{\lambda}}  (A_1\otimes A_2 \vert B_1 \otimes B_2) \geq \lambda \widetilde{T_{\lambda}}  (A_1 \vert B_1)\otimes \widetilde{T_{\lambda}}  ( A_2 \vert B_2).  $
\item[(b')] $\widetilde{T_{\lambda}}  (A_1\otimes A_2 \vert B_1 \otimes B_2) \leq  \widetilde{T_{\lambda}}  (A_1 \vert B_1)\otimes A_2 + A_1 \otimes \widetilde{T_{\lambda}}  ( A_2 \vert B_2).$
\end{itemize}
\item[(2)]  For any $\lambda \leq 0$ and $0< B_i \leq A_i$, $(i=1,2)$, we have the following inequalities.
\begin{itemize}
\item[(c')] $\widetilde{T_{\lambda}}  (A_1\otimes A_2 \vert B_1 \otimes B_2) \leq \lambda \widetilde{T_{\lambda}}  (A_1 \vert B_1)\otimes \widetilde{T_{\lambda}}  ( A_2 \vert B_2).  $
\item[(d')] $\widetilde{T_{\lambda}}  (A_1\otimes A_2 \vert B_1 \otimes B_2) \leq  \widetilde{T_{\lambda}}  (A_1 \vert B_1)\otimes A_2 + A_1 \otimes \widetilde{T_{\lambda}}  ( A_2 \vert B_2).$
\end{itemize}
\end{itemize}
\end{Prop}

\begin{Rem}
(a') and (c') are the same inequalities with (a) and (c) in Corollary \ref{fininal_cor}. 
We easily find from the above (b') and (d') that we have the subadditivity without depending on the ordering $A_i$  and $B_i$. 
By the similar way of Remark \ref{final_rem}, this subadditivity for Tsallis relative operator entropy implies the subadditivity for Tsallis entropy:
\begin{equation}  \label{rem2_eq}
S_{\lambda}(\rho_1 \otimes \rho_2) \leq S_{\lambda}(\rho_1) + S_{\lambda}( \rho_2),
\end{equation}
for any $\lambda \leq 0$ and any invertible density operators $\rho_i,(i=1,2)$. 
\end{Rem}

In Tsallis statistical physics, the subadditivity Eq.(\ref{rem2_eq}) and the superadditivity Eq.(\ref{cor2}) are famous and fundamental.
Considering the tensor product, we could show the subadditivity ((b') and (d') in Proposition \ref{final_prop}) 
and the superadditivity ((b) and (d) in Corollary \ref{fininal_cor}) for Tsallis relative operator entropy as an operator inequality version, respectively.
Our results can be seen as the generalizations of Eq.(\ref{rem2_eq}) and Eq.(\ref{cor2}).


\begin{thebibliography}{99}
\bibitem{Ab:re1} S.Abe, Monotonic decrease of the quantum nonadditive divergence by projective measurements, 
Physics Letters A, {\bf 312}(2003),336--338.

\bibitem{An} T.Ando, {\em Topics on operator inequality}, Lecture Notes, Hokkaido Univ.,Sapporo,1978.

\bibitem{F} J.I.Fujii, Operator means and the relative operator entropy,
Operator Theory: Advances and Applications, {\bf 59} (1992), 161--172.

\bibitem{FFS} J.I.Fujii, M.Fujii and Y.Seo, An extension of the Kubo-Ando theory: Solidarities, 
 Math.Japonica, {\bf 35} (1990), 387--396.

\bibitem{FuKa:rel} J.I.Fujii and E.Kamei, Relative operator entropy in noncommutative information theory, 
 Math. Japonica, {\bf 34}(1989)341-348.

\bibitem{FK} J.I.Fujii and E.Kamei, Uhlmann's interpolational method for operator means, 
 Math.Japonica, {\bf 34}(1989), 541-547.

\bibitem{FYK} S.Furuichi, K.Yanagi and K.Kuriyama, Fundamental properties for Tsallis relative entropy, 
J.Math.Phys., {\bf 45}(2004),4868--4877.

\bibitem{YKF} K.Yanagi, K.Kuriyama and S.Furuichi, Generalized Shannon inequalities based on Tsallis relative operator entropy, 
 Linear Algebra Appl., {\bf 394}(2005), 109--118. 

\bibitem{Furuta} T.Furuta, {\em Invitation to Linear Operators: From Matrix to Bounded Linear Operators on a Hilbert Space}, CRC Pr I Llc, 2002.

\bibitem{HK} F.Hiai and H.Kosaki, {\em Means of Hilbert Space Operators}, Springer, 2003.

\bibitem{HP} F.Hiai and D.Petz,  The proper formula for relative entropy in asymptotics in quantum probability, 
Comm.Math.Phys.,{\bf 143}(1991), 99--114.

\bibitem{NU} M.Nakamura and H.Umegaki, A note on the entropy for operator algebras,  Proc. Jap. Acad., {\bf 37} (1961),149--154.

\bibitem{Re} A.R\'enyi, On the foundation of information theory,  Rev.Int.Stat.Inst., {\bf 33}(1965),1--14.

\bibitem{Ts:rel} C.Tsallis, Possible generalization of Boltzman-Gibbs statistics,  J.Stat. Phys., {\bf 52}(1988),479--487.

\bibitem{AO:re1} C.Tsallis et al., {\em Nonextensive Statistical Mechanics and Its Applications}, edited by S. Abe and Y. Okamoto (Springer-Verlag, Heidelberg,2001); 
see also the comprehensive list of references at \texttt{http://tsallis.cat.cbpf.br/biblio.htm}

\bibitem{Ume} H.Umegaki,  Conditional expectation in an operator algebra, IV (entropy and information),
Kodai Math.Sem.Rep., {\bf 14}(1962), 59--85.

\bibitem{Zh} X.Zhan, {\em Matrix Inequalities}, Springer, 2002.

\end{thebibliography}
\end{document}